\documentstyle[twoside,12pt]{article}
\input psfig
\def\IMSmarkvadjust{0 pt}
\def\IMSmarkhadjust{0 pt}
\def\IMSmarkhpadding{0 pt}
\def\IMSpubltext{Published in modified form:}
\def\SBIMSMark#1#2#3{
 \font\SBF=cmss10 at 10 true pt
 \font\SBI=cmssi10 at 10 true pt
 \setbox0=\hbox{\SBF \hbox to \IMSmarkhpadding{\relax}
                Stony Brook IMS Preprint \##1}
 \setbox2=\hbox to \wd0{\hfil \SBI #2}
 \setbox4=\hbox to \wd0{\hfil \SBI #3}
 \setbox6=\hbox to \wd0{\hss
             \vbox{\hsize=\wd0 \parskip=0pt \baselineskip=10 true pt
                   \copy0 \break%
                   \copy2 \break%
                   \copy4 \break}}
 \dimen0=\ht6   \advance\dimen0 by \vsize \advance\dimen0 by 8 true pt
                \advance\dimen0 by -\pagetotal
	        \advance\dimen0 by \IMSmarkvadjust
 \dimen2=\hsize \advance\dimen2 by .25 true in
	        \advance\dimen2 by \IMSmarkhadjust

  \openin2=publishd.tex
  \ifeof2\setbox0=\hbox to 0pt{}
  \else 
     \setbox0=\hbox to 3.1 true in{
                \vbox to \ht6{\hsize=3 true in \parskip=0pt  \noindent  
                {\SBI \IMSpubltext}\hfil\break
                {\it  Math. Intelligencer}~{\bf 17} (1995), 48--56
 
                \vfill}}
  \fi
  \closein2
  \ht0=0pt \dp0=0pt
 \ht6=0pt \dp6=0pt
 \setbox8=\vbox to \dimen0{\vfill \hbox to \dimen2{\copy0 \hss \copy6}}
 \ht8=0pt \dp8=0pt \wd8=0pt
 \copy8
 \message{*** Stony Brook IMS Preprint #1, #2. #3 ***}
}

\newtheorem{lemma}{Lemma}

\begin{document}
\title{\centerline{\psfig{figure=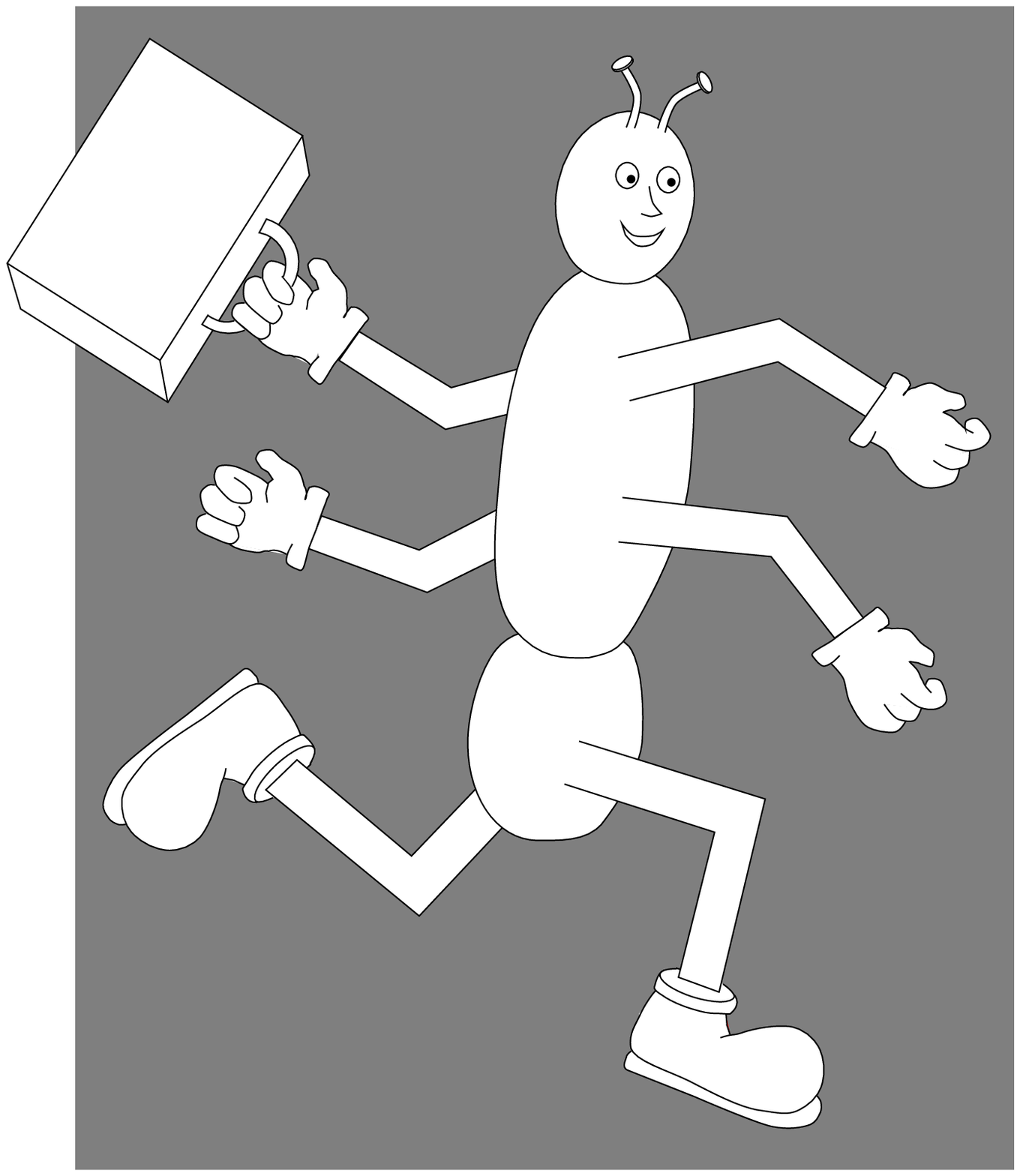,height=.75in}}
       Further Travels with My Ant}
\author{David Gale, Jim Propp, Scott Sutherland, and Serge Troubetzkoy}
\date{}
\maketitle
\SBIMSMark{1995/1}{January 1995}{}

\thispagestyle{empty}


\section{Introduction}

   A recurring theme of this column over the past four years has been
what we have referred to as computer-generated mysteries.  Examples are
sequences defined by simple rational recursions whose terms turn out
to be integers with interesting but unexplained divisibility
properties, or geometric configurations that are observed to exist
although there are no proofs of existence.  In most of the examples
the reported mysteries have remained unsolved and in some cases may in
fact be, in a suitable sense, ``unsolvable''.  It is therefore gratifying
to be able to present in this month's column an
elegant solution of a previously described mystery.  An especially
pleasing feature of this solution is the fact that the ``breakthrough''
was made possible by drawing the right picture.  Once the picture is
drawn it becomes clear what must be proved, after which further study
of the pictures gives the clue as to how to construct the proof.  It
turns out that at one point one needs to use the Jordan curve theorem
for a special class of closed curves.  

    In the paragraphs to follow we shall present an informal but
hopefully convincing argument which will require little more of
the reader than the careful observation of a collection of pictures.

\section{The Story So Far and the Mystery}

   We shall be concerned once again with a certain automaton that has
been called an ant.  An ant lives in a plane that has been partitioned
in the standard way into squares, or, as we shall call them, {\it cells}
(think of the corners of the cells as the lattice points of the
plane).  Each cell can be in one of several states, and the respective
states of these cells change over the course of time, as the cells are
in turn visited by the ant.  Before we explain precisely how the ant
interacts with its surroundings, it will help to consider a specific
example.  In the simplest non-trivial version of the ant, originally
studied by Chris Langton, there are just two states, which we will
call L (for left) and R (for right).  
We think of the ant as starting on the boundary between two cells, heading 
in one of the four cardinal compass directions (East or West on the vertical
boundaries, and North or South on the horizontal ones). As the ant advances
through the cell, it makes a 90~degree turn, turning left in L-cells and
right in R-cells, and moves to the boundary of the the neighboring cell.  As
the ant leaves each cell, it changes the state of the cell, switching
L-cells to R-cells, and vice versa. 
For an informal discussion of this ant see [1], where a number of interesting 
ant behaviors are observed.  For our present purposes we note only the
following still unexplained phenomenon: If initially all cells were in
the same state, then at
various times the ant's ``track'' (the set of cells it has visited, together
with their current states) is centrally symmetric, that is, the configuration
of L- and R-cells has central symmetry.  We reproduce here, in Figure~1, the
corresponding pictures, where the R-cells are black and the 
L-cells white; all cells were initially in the L state and the ant started
out heading west on the right boundary of the central cell, and about to
enter it.
These central symmetries stop occurring eventually
and after about 10,000 time-units the ant settles into a periodic
``highway-building'' behavior, heading off to infinity in a
southwesterly direction.  This phenomenon of ``transient symmetry''
awaits a satisfying explanation.

\begin{figure}[htp]
\centerline{\psfig{figure=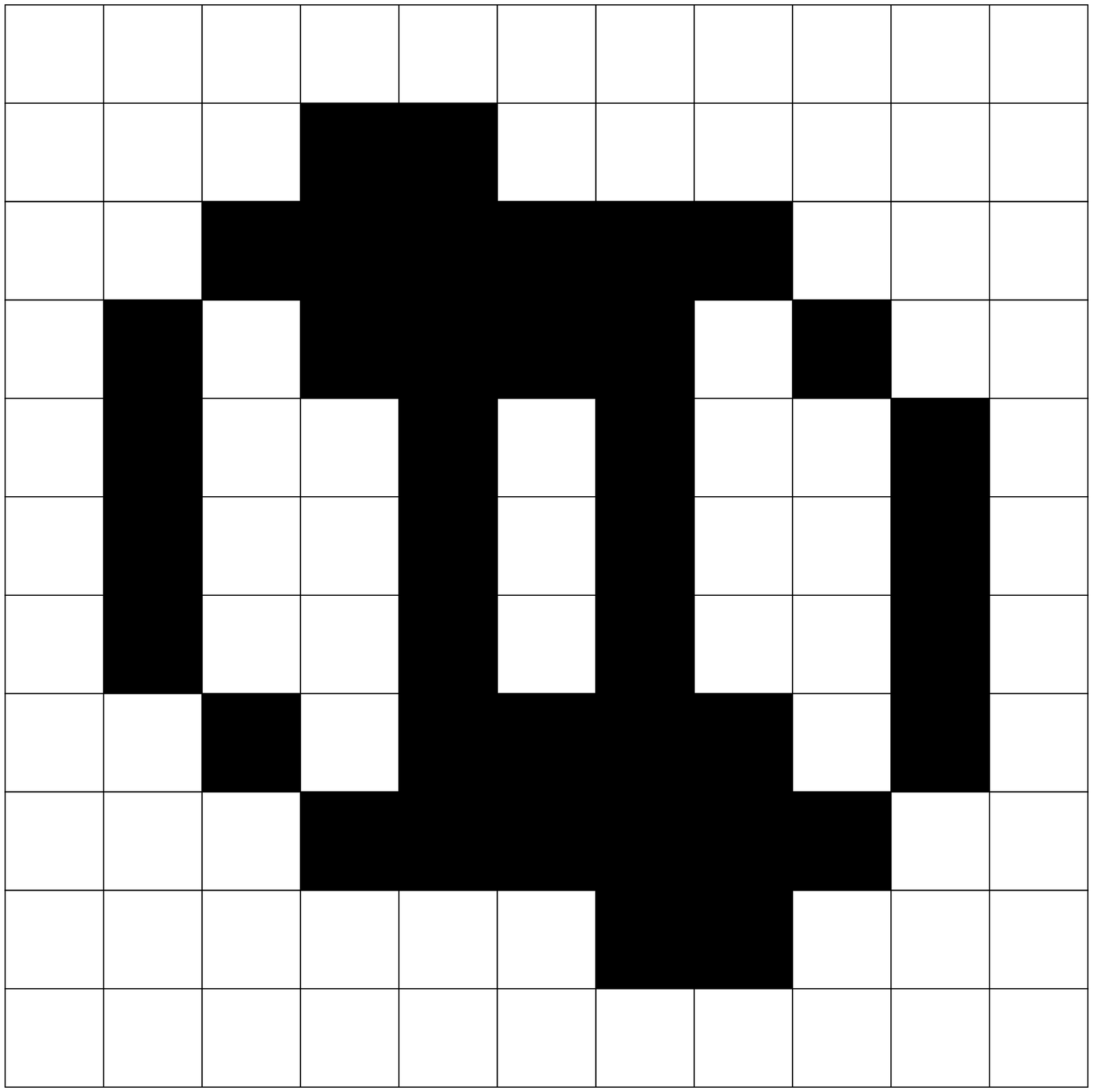,height=1.5in}\hfil
	    \psfig{figure=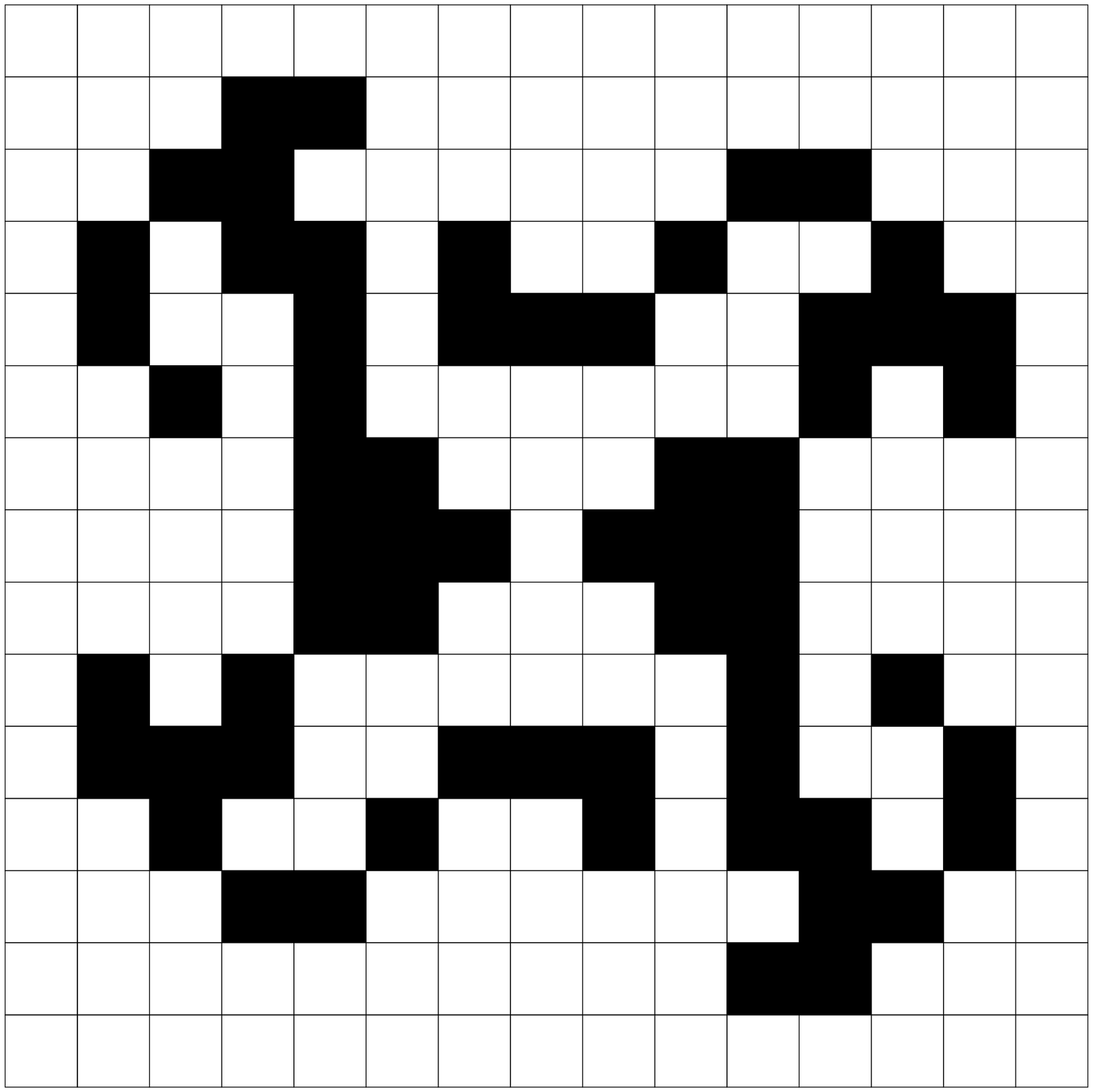,height=1.5in}\hfil
	    \psfig{figure=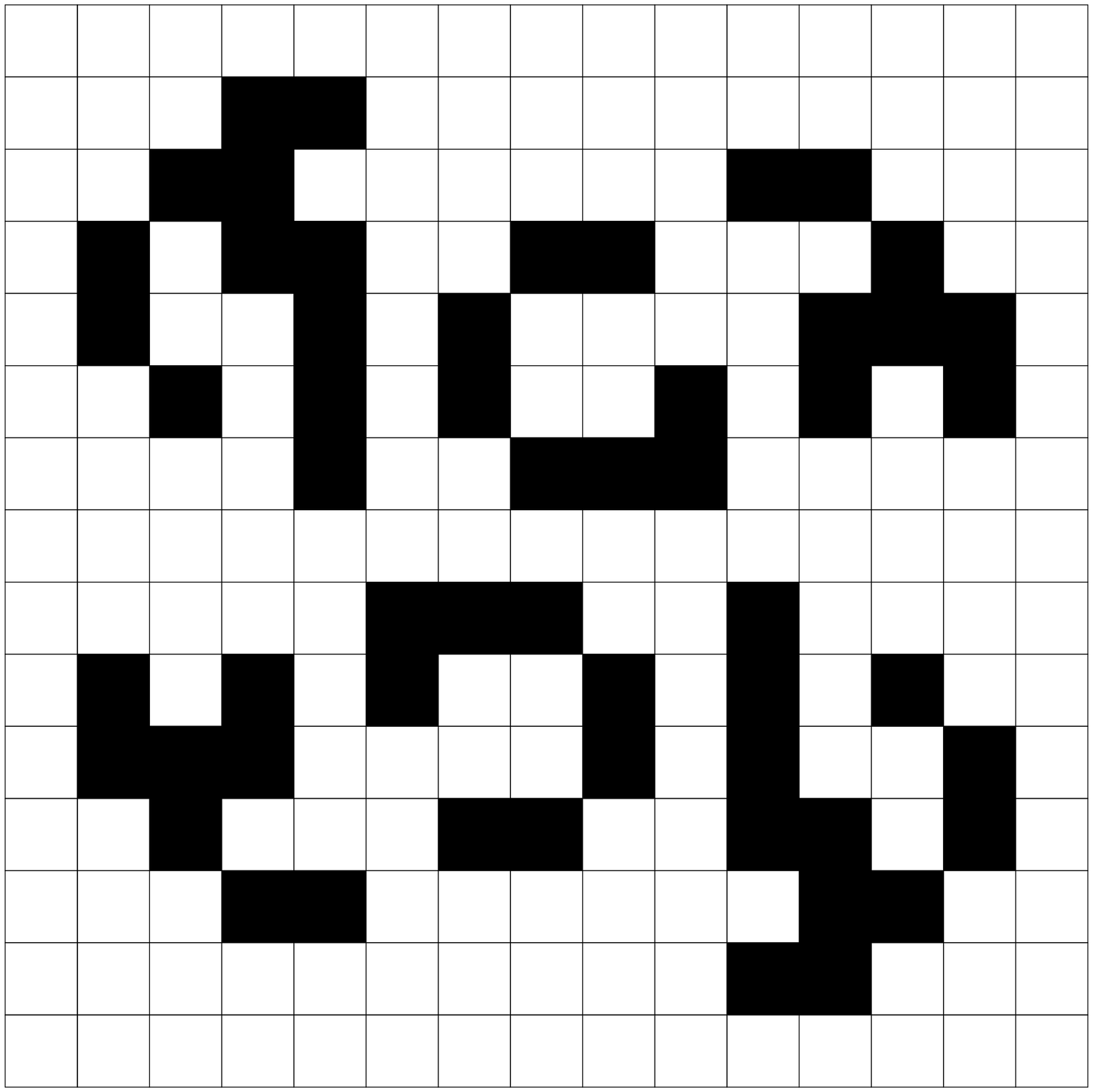,height=1.5in}}
\caption{The universe of ant 2 at times 184, 368 and 472.}
\end{figure}

In our next episode [2], Jim Propp describes the behavior of the more 
general $n$-state ants.  There are now $n$ different states for cells to be in,
numbered 1 through $n$, and it is the ant's internal ``programming'' which
tells it which states to treat as L's and which states to treat as R's
in choosing which way to turn.  We represent this programming by means
of a {\it rule-string}, a string of $n$ L's and R's whose $k$-th letter
represents the ant's course of action when it comes to a cell in state
$k$; for instance, the 7-state string LLRRRLR represents an ant that
turns left on visiting cells in states 1, 2, and 6 and turns right
on visiting cells in states 3, 4, 5, and 7.  When an ant comes to an
L-cell (a cell in an L-state), it turns left; when it comes to an
R-cell, it turns right.  When the ant leaves a cell that is in state $i$
the cell changes to state $i+1$ (mod $n$).  In this terminology the simple
ant has the rule-string LR.

{}From symmetry considerations it clearly suffices to consider
only strings that start with an L.  If in the rule-string we replace
an L by a 1 and an R by 0 we see that each ant corresponds to a
positive integer expressed in base two, so the simple ant is ant 2 and
our seven-state ant with ``genome'' LLRRRLR is ant 98.  Propp finds that
different ants behave very differently depending on their
rule-strings.  Some seem to be completely chaotic, while others
eventually build highways.  The only general statement that can be
made is the following

\bigskip
\begin{center}
{\bf Fundamental Theorem of Myrmecology}
(Bunimovich-Troubetzkoy):\\
{\it  An ant's track is always unbounded.}
\end{center}

\bigskip  (For a proof see [1].  Note that the attribution given there 
is incorrect; the proof actually first appeared in [3].)
 
    In this exposition we concentrate on the phenomenon described by
Propp as follows.

     ``Ants 9 and 12 [LRRL and LLRR in our notation] are the truly
surprising ones [among ants with rule-strings of length 4].  In each
case the the patterns get ever larger, but without ever getting too
far away from bilateral symmetry!  More specifically, one finds that
the ant makes frequent visits to the cell it started from, and when it
does, the total configuration quite often has bilateral symmetry.''

\begin{figure}[htp]
\hbox to \hsize{%
  \vbox{\hsize=.48\hsize
	\centerline{\psfig{figure=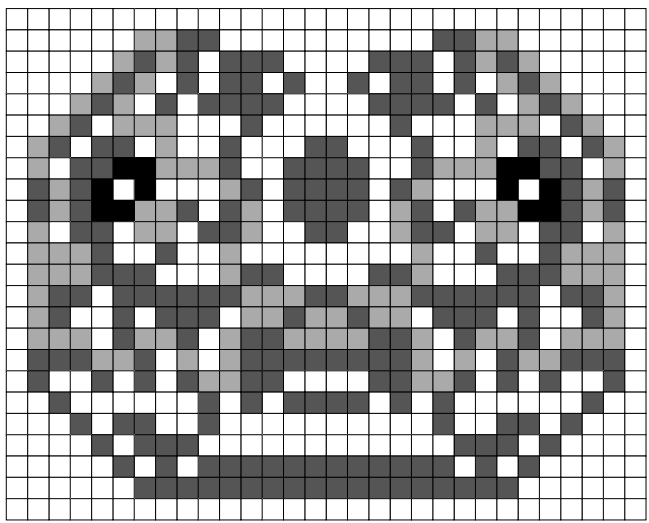,height=\hsize}}
	\caption{The universe of ant 12 at time 16,464.}}
  \hfil     
  \vbox{\hsize=.48\hsize
	\centerline{\psfig{figure=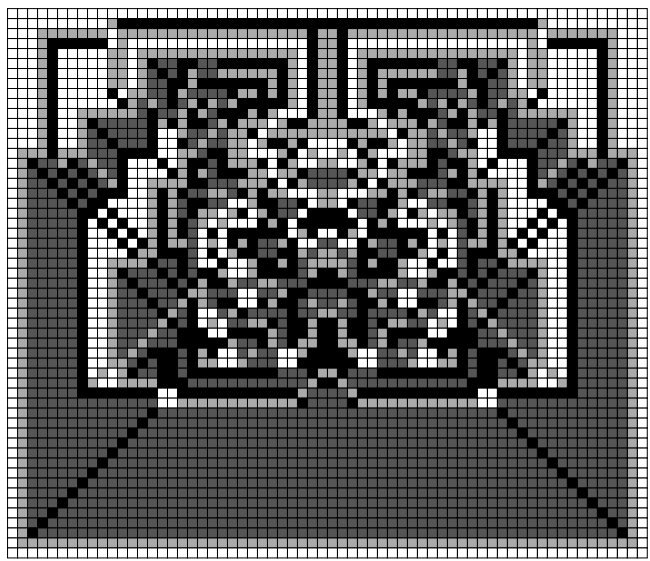,height=\hsize}}
	\caption{The universe of ant 9 at time 38,836. }}
}
\end{figure}

Figures 2 and 3, reproduced from [2], show the track of ant 12 after
16,464 steps and the track of ant 9 after
36,836 steps.  The figures use black, white and two shades of grey to
indicate the four different states.  Propp continues, ``To find more
ants of this sort we have to move to rule strings of length 6.  Here
we encounter another mystery: the rule strings that lead to
bilaterally symmetric patterns are 33, 39, 48, 51, and 60.  Note that
all these numbers are divisible by three!  Surely this cannot be an
accident.''  Indeed it is not, as we will presently show.  
 
\section{The Picture. Truchet Tiles.}

\begin{figure}[htp]
	\centerline{\psfig{figure=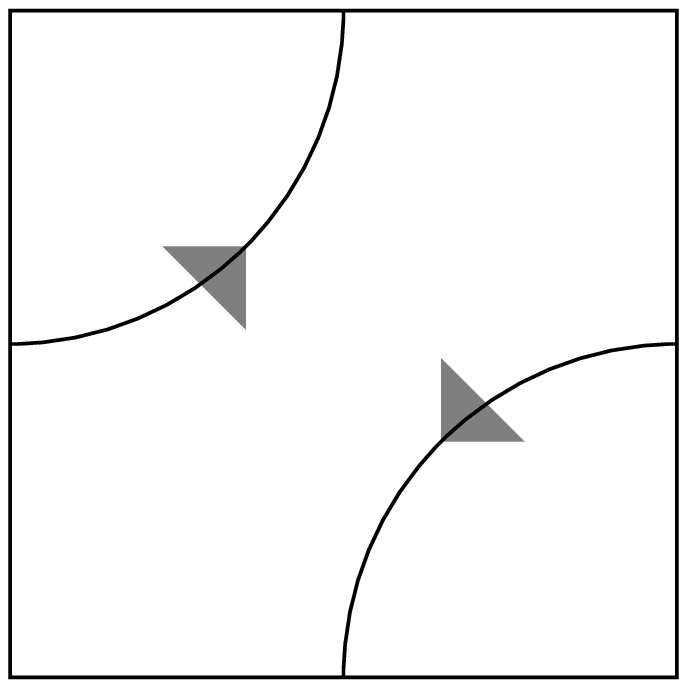,height=1in} \hskip .5 in
		    \psfig{figure=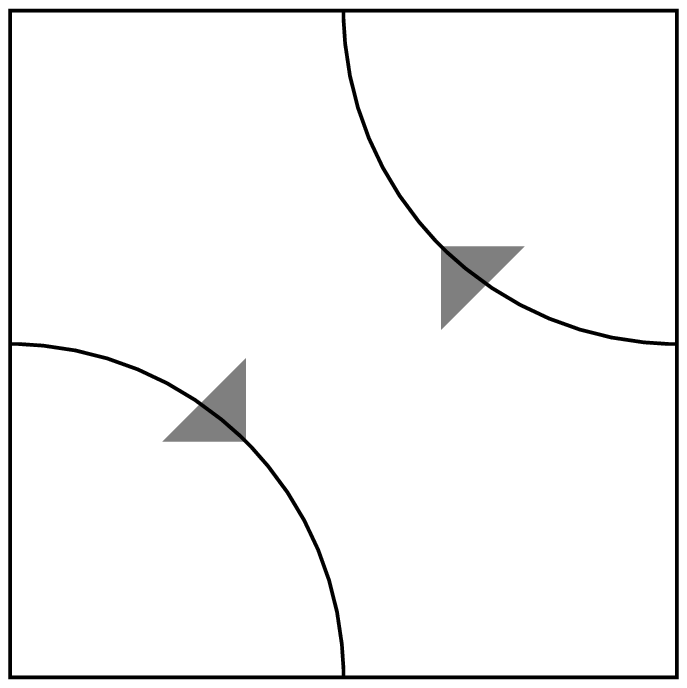,height=1in}}
	\caption{ H-cells in L and R orientation. }
\bigskip
	\centerline{\psfig{figure=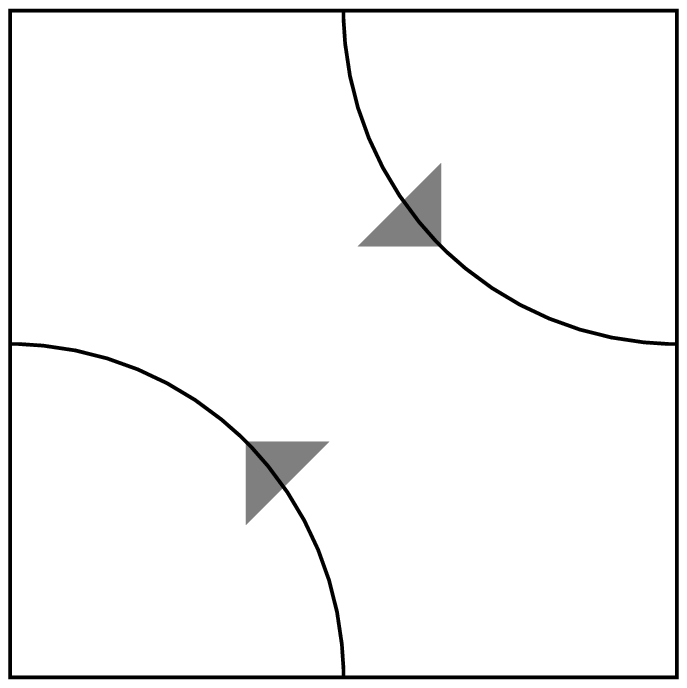,height=1in}\hskip .5 in
		    \psfig{figure=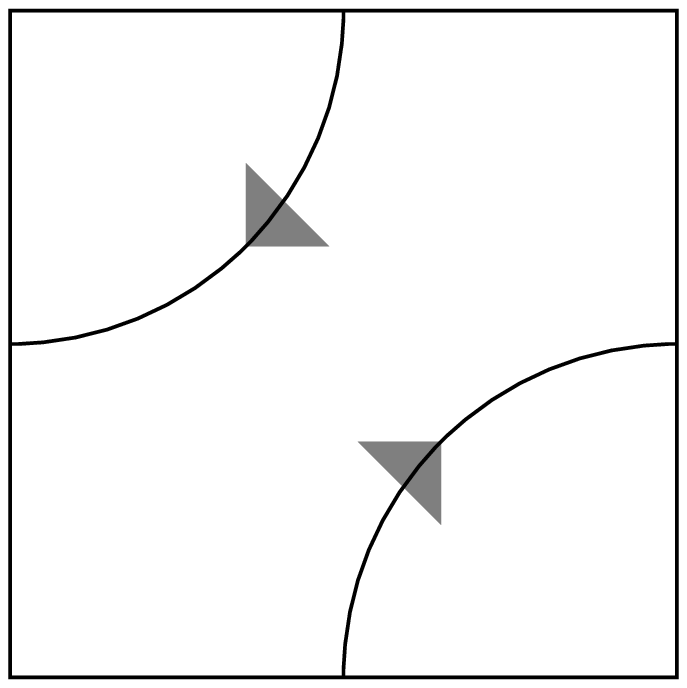,height=1in}}
	\caption{V-cells in L and R orientation. }
\end{figure}

   First a preliminary observation.  Since an ant turns through a
right angle after each move, it follows that its moves are alternately
horizontal and vertical.  Therefore the cells of the plane split up
checkerboard-fashion into two sets: H-cells, which the ant always
enters  horizontally, from the left or the right, and exits
vertically, either at the top or the bottom; and V-cells, which the
ant always enters vertically and exits horizontally.  Now, a crucial
property of a cell is that it is able to change its state and thus
change the course of action (left-turn versus right-turn) that the ant
will take after its next visit to the cell.  The key idea for
representing this property, due to Bernd R\"ummler, was to actually
install schematic ``switches'' in the cells, in which the orientation of
the switch indicates whether the cell is currently an L-state or an
R-state.  These switches are called Truchet tiles and are illustrated
in Figure 4, which shows an H-tile first in its L and then in its R
orientation, and Figure 5 which shows the same thing for a V-tile.
(An ``H-tile'' is the Truchet tile associated with an H-cell, and similarly
for V-tiles.)  Note that switching orientation corresponds to
reflecting the tile about either its vertical or horizontal axis.

\begin{figure}[ht]
	\centerline{\psfig{figure=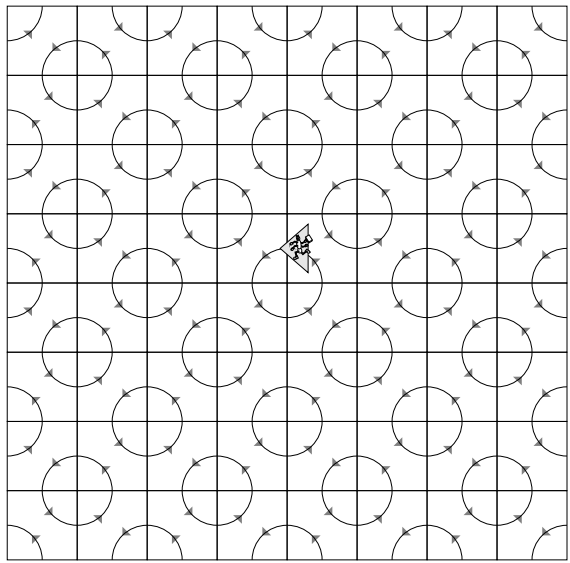,width=.42\hsize}}
	\caption{The initial configuration.}
\end{figure}

Figure 6 shows the plane paved with Truchet tiles all of which
have the L orientation, giving a pattern of disjoint circles: the
``initial configuration''.  As the ant moves, some of the tiles will
switch in accordance with the given rule-string, but it is clear that
the pattern will always be made up of a set of disjoint simple closed
curves which we shall call {\it Truchet contours}.  (Infinite contours
cannot arise, since at each stage the ant has only switched finitely
many of the tiles.)  It will be helpful (at least for now) to imagine
that the ant actually travels along the Truchet curves themselves,
rather than along a lattice-path joining the centers of the cells, and
that the ant's initial position is the midpoint of one of the edges of
the central Truchet tile.  Figure 7 gives the ``Truchet pictures''
corresponding to the (transiently) centrally-symmetric configurations
of the simple ant given by the earlier Figure 1, while Figures 8 and 9
give Truchet pictures corresponding to the earlier Figures 2 and 3,
and manifest the same bilateral symmetry.  In all of these pictures 
the ant has returned to its original location, and of particular interest
is the Truchet contour through this point, which we will call the
{\it principal contour}.  Initially, the principal contour is just a circle
as are all other contours. In Figure 9 the principal contour has been
highlighted. 

\begin{figure}[htp]
	\centerline{\psfig{figure=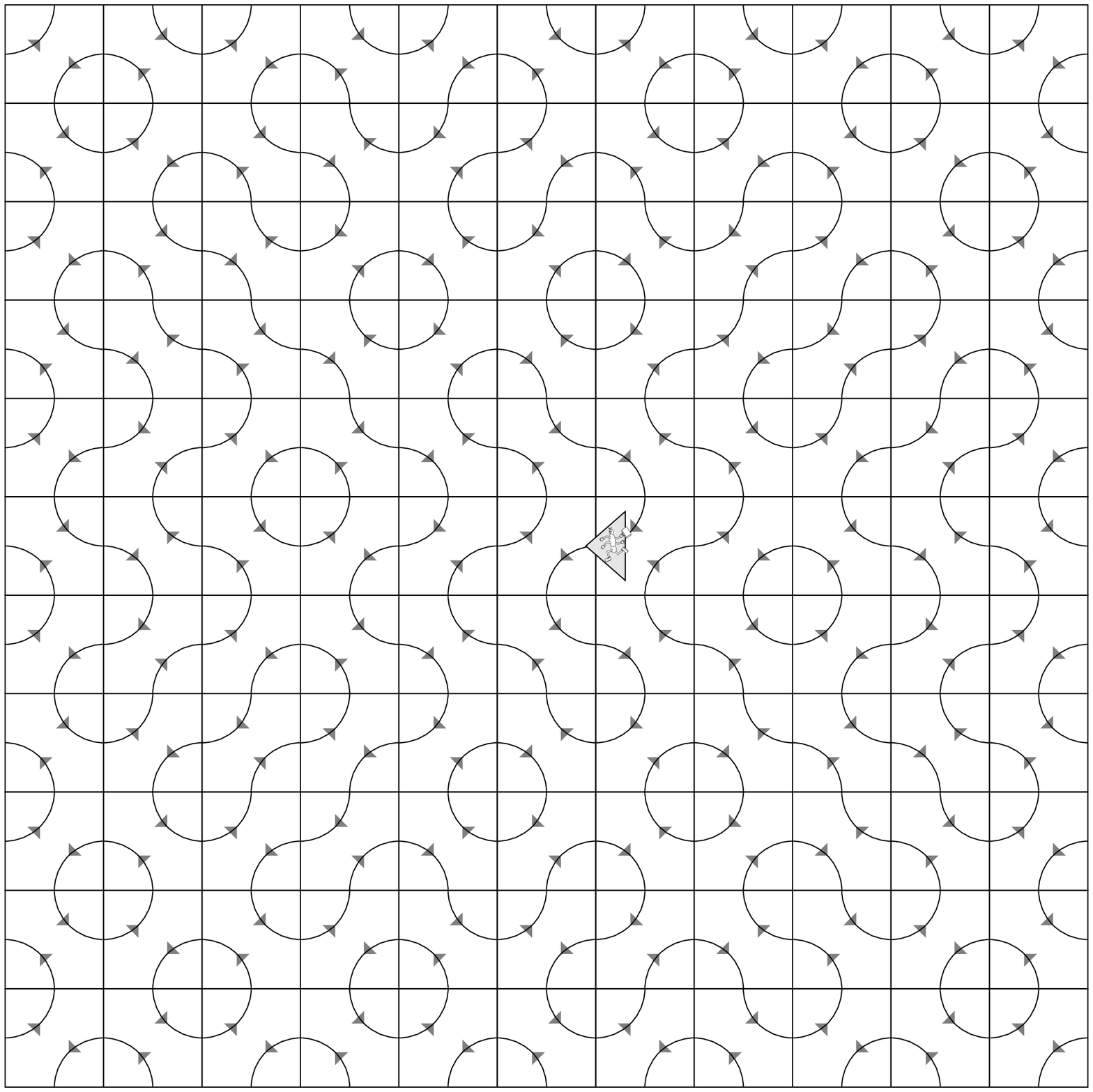,width=.45\hsize}}
	\bigskip
	\centerline{\psfig{figure=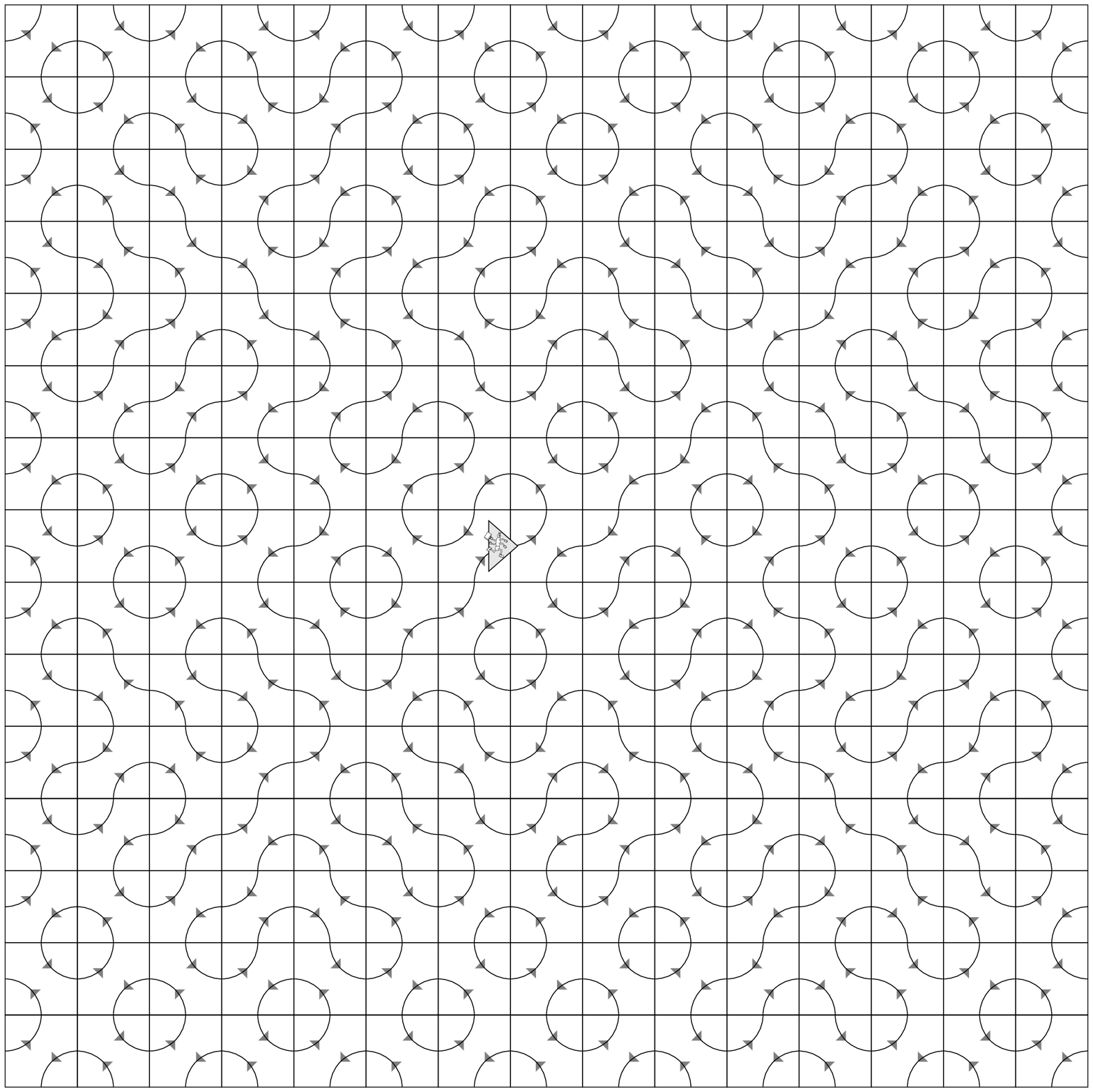,width=.45\hsize} \hfil
	            \psfig{figure=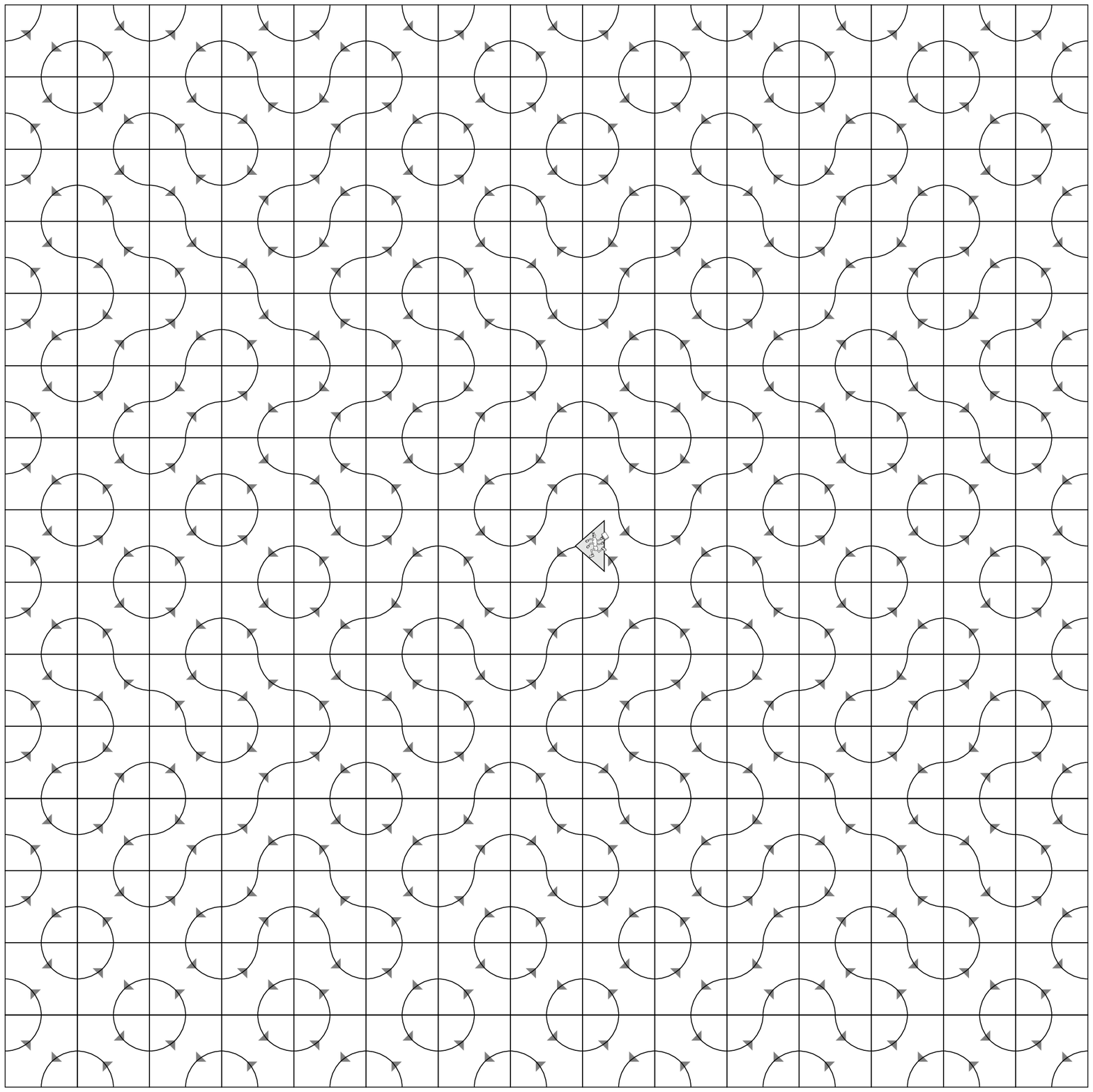,width=.45\hsize}}
	\caption{ The universe of ant 2 at times 184, 368 and 472. }
\end{figure}

\begin{figure}[htp]
	\centerline{\psfig{figure=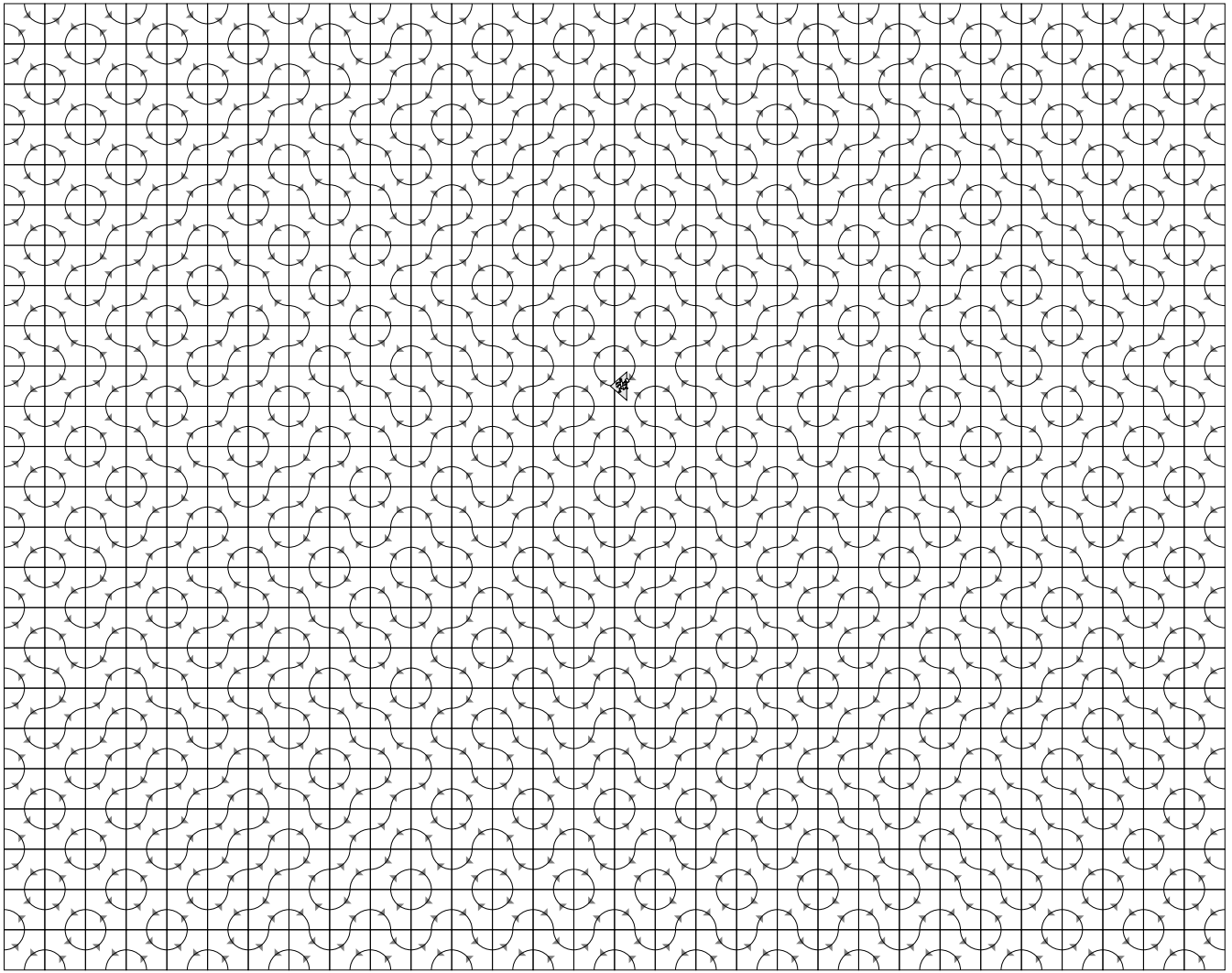,width=\hsize}}
	\caption{ The universe of ant 12 at time 16,464 }
\end{figure}

\begin{figure}[htp]
	\centerline{\psfig{figure=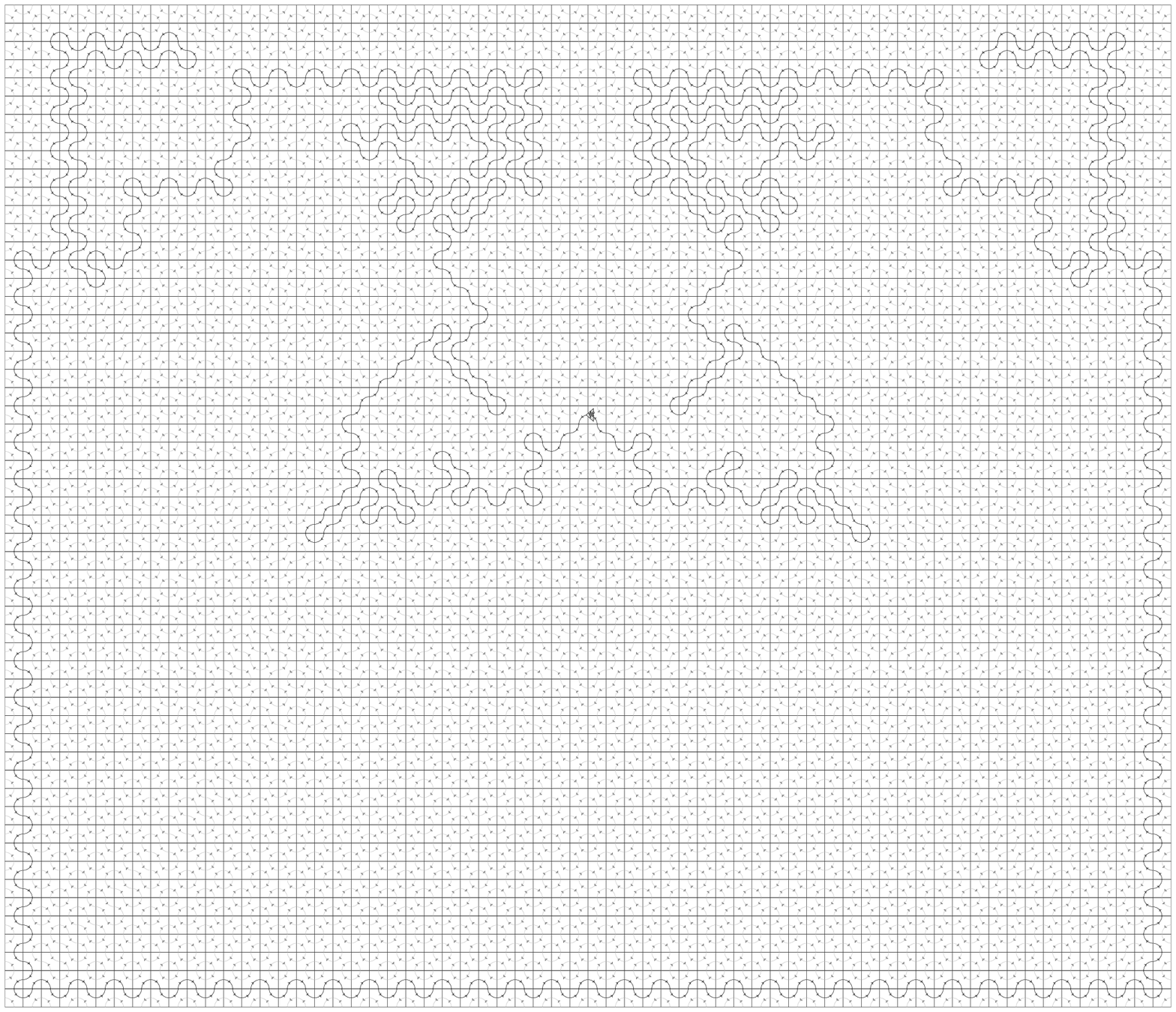,width=\hsize}}
	\caption{The universe of ant 9 at time 38,836. The principal
                 contour has been highlighted.} 
\end{figure}

Now if one knew that each time the ant left its starting location
it would stay on the principal contour until it returned once more to
its starting location then it would follow that when the ant had
completed its tour all the symmetries of the initial state of the
universe would be preserved; for whenever a cell was visited its
symmetric mate would also be visited. 

\begin{figure}[htp]
	\centerline{\psfig{figure=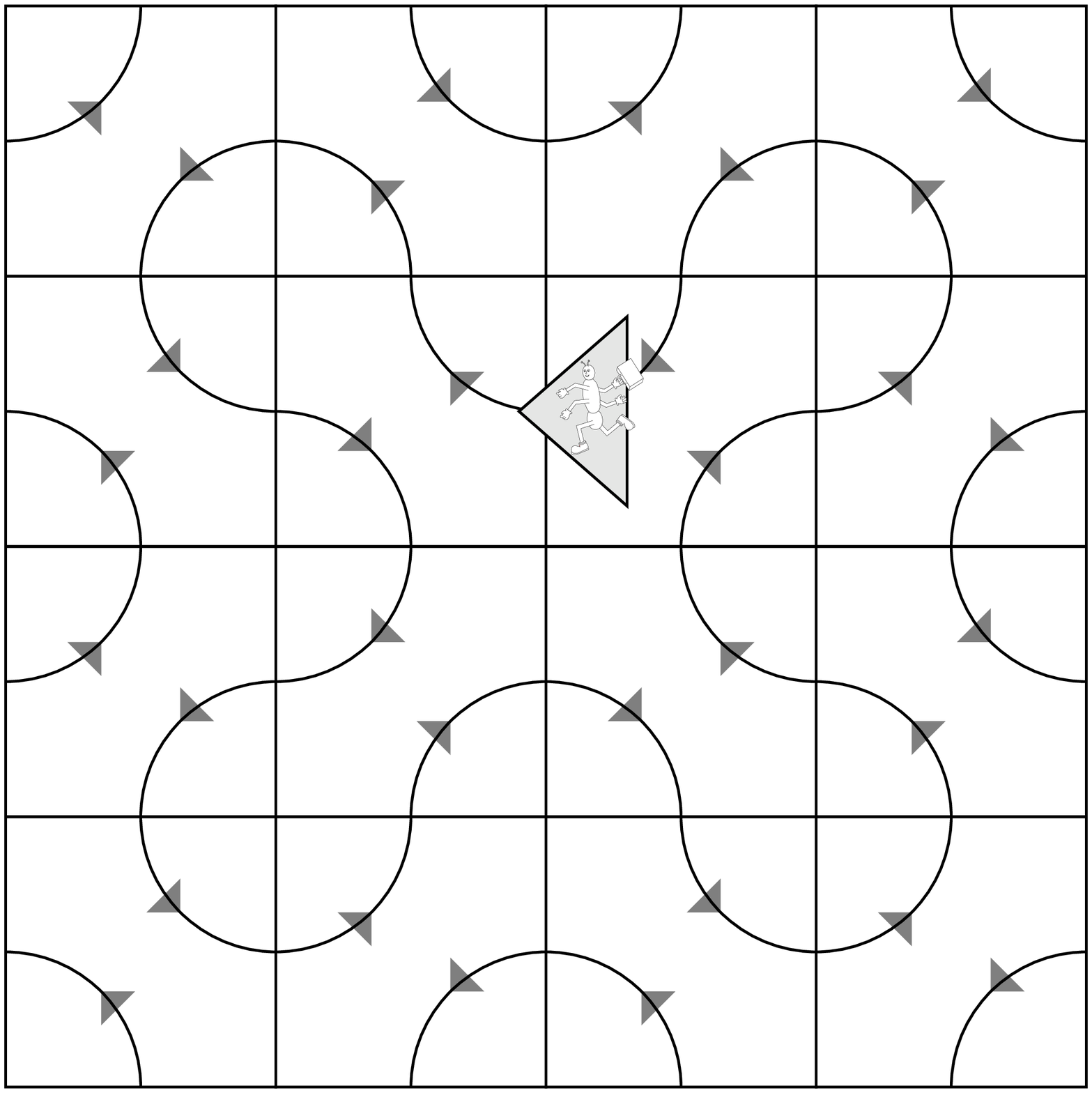,width=1.5in}\hfil
		    \psfig{figure=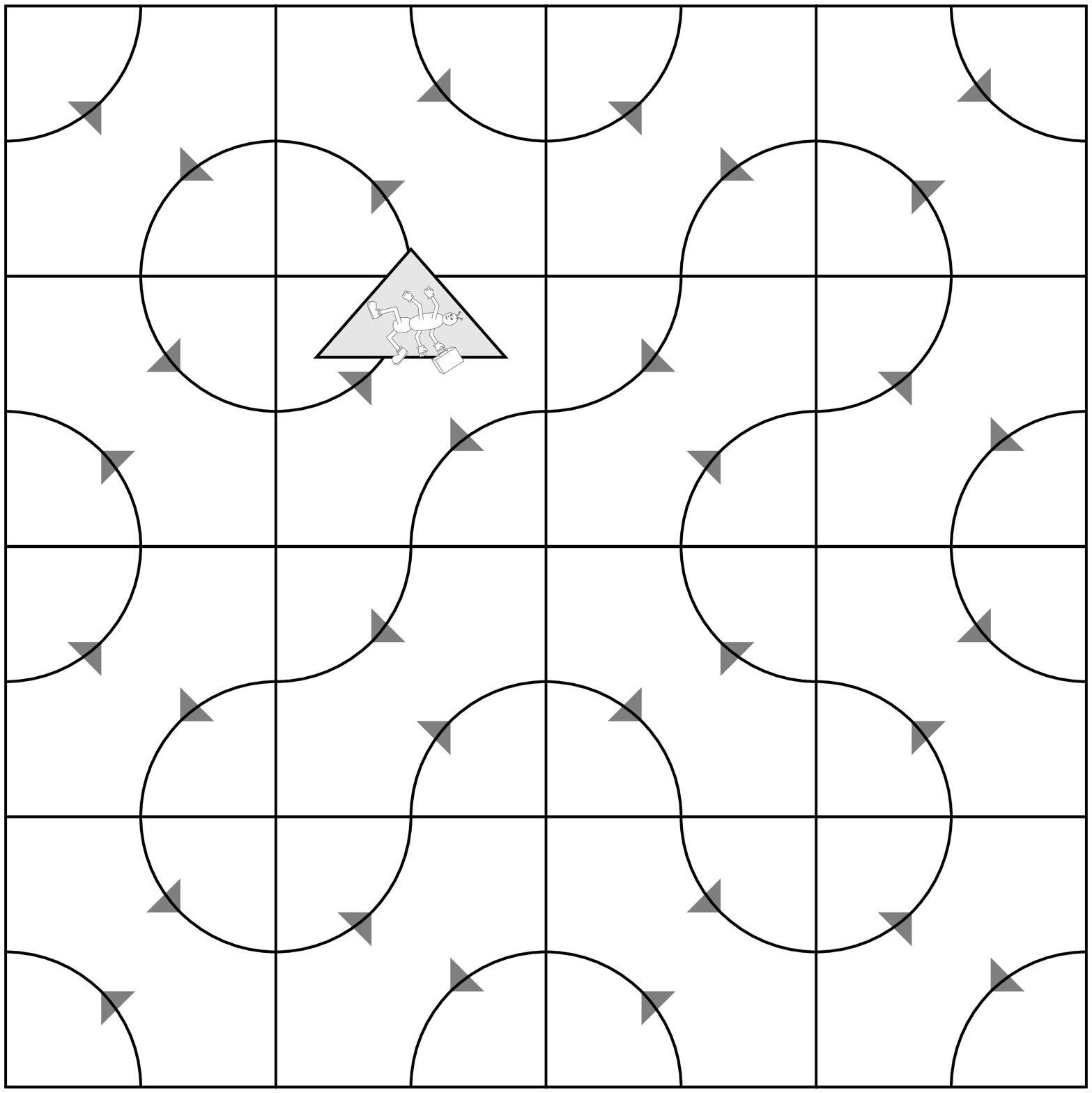,width=1.5in}}
	\caption{Switching: the universe of ant 2 at times 4 and 5.}
\end{figure}

   Unfortunately, an ant will not in general stay on the principal
contour, because (as in Figure 10) the contour may pass through some cells
twice.  This means that when the ant returns to such a cell the
orientation of that cell may have switched (indeed this will always be
the case for the simple ant studied by Langton), causing the ant to
leave the principal contour.  However, for the ants listed by Propp
one can show that in fact

\begin{equation}
   {\mbox{\it The cells that are visited twice by a contour never
   switch on the first visit. }}
\end{equation}   

    This implies that the ant will engage in a process of repeatedly
tracing out bilaterally-symmetric principal contours, resulting in a
bilaterally-symmetric universe each time the ant returns to its
starting point.  Property (1) was first proved by R\"ummler for ant 12
and then generalized to the other ants by Troubetzkoy.  The argument
to follow is a reworking by Propp of those proofs.

\section{The Even Run-Length Property and the Augmented Picture}

   Why do some ant tracks exhibit recurrent bilateral symmetry and
others not?  Here are the rule-strings for the 4- and 6-state ants of
Propp's article that exhibit recurrent symmetry. 

\begin{center}
\begin{tabular}{ll}\\
{\bf Ant}  &  {\bf Rule-string}\\
 9  &   LRRL\\
12  &  LLRR \\
33  & LRRRRL\\
39  &   LRRLLL\\
48  &   LLRRRR\\
51  &  LLRRLL\\
57  & LLLRRL\\
60  & LLLLRR
\end{tabular}
\end{center}

  It is not hard to see what these strings have in common.  If we
think of them arranged in cyclic rather than linear order then each of
them consists of an even number of L's followed by an even number of
R's.  In general we say a rule-string has the 
{\it even run-length property}
if in the cyclic order it consists of alternate runs of L's and R's of
even length, e.g.\ LRRLLRRRRL. 
For simplicity in what follows we will consider only the
case where the string starts with an even number of L's (ants 12, 48,
51 and 60), the argument for the other case being similar.  

\begin{figure}[htp]
	\centerline{\psfig{figure=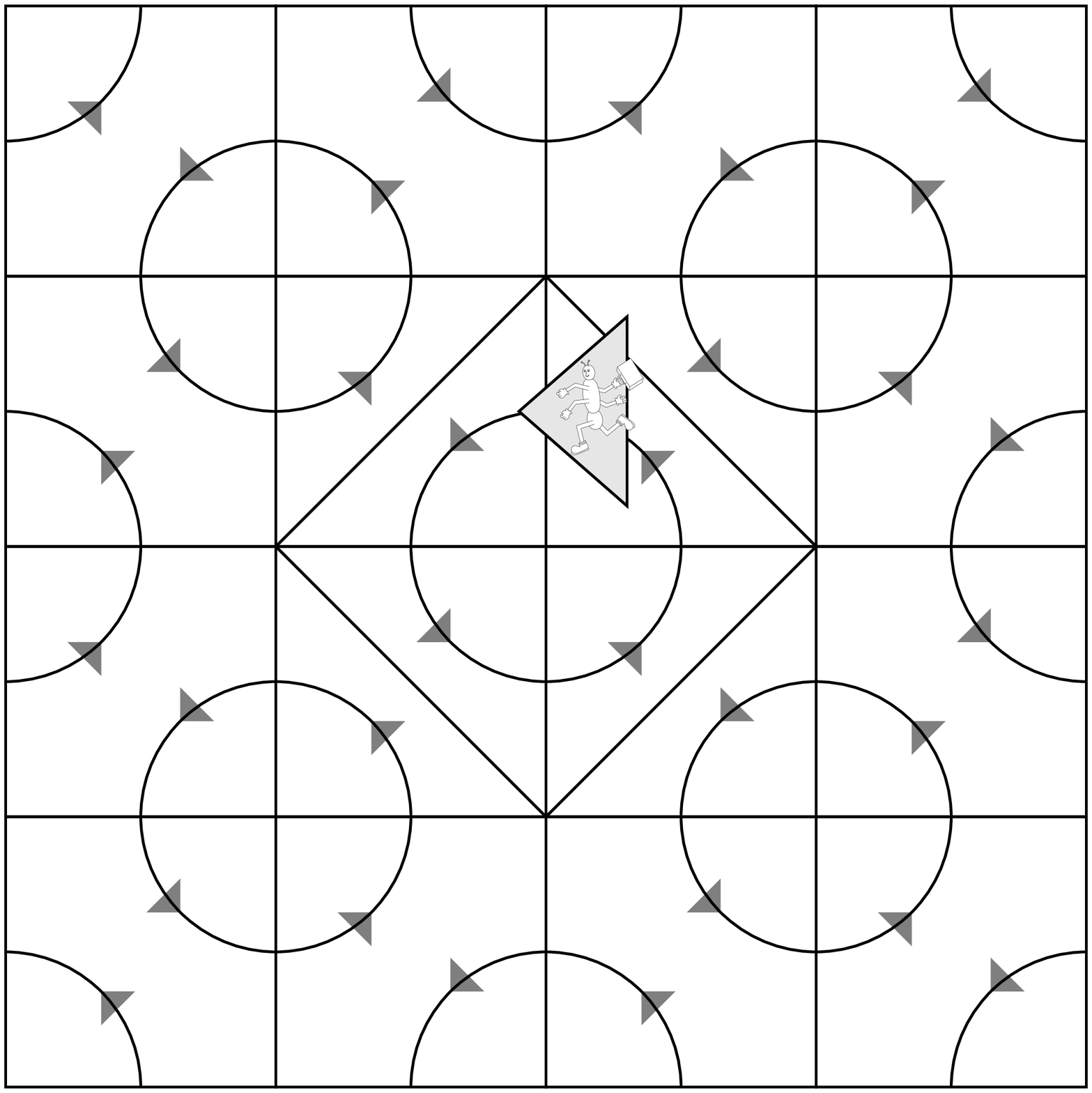,height=2.5in}\hfil
		    \psfig{figure=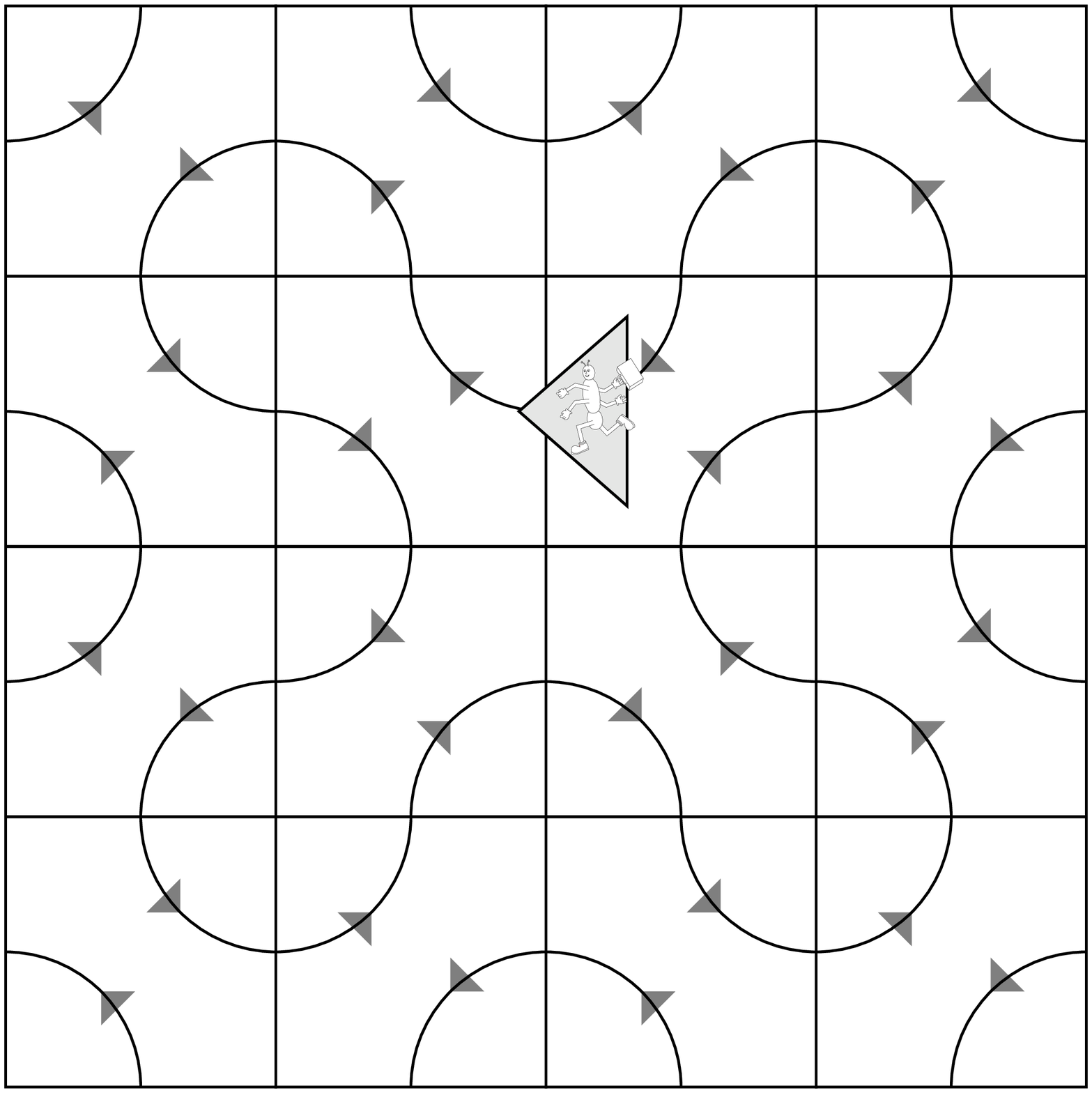,height=2.5in}}
	\bigskip\bigskip\bigskip
	\centerline{\psfig{figure=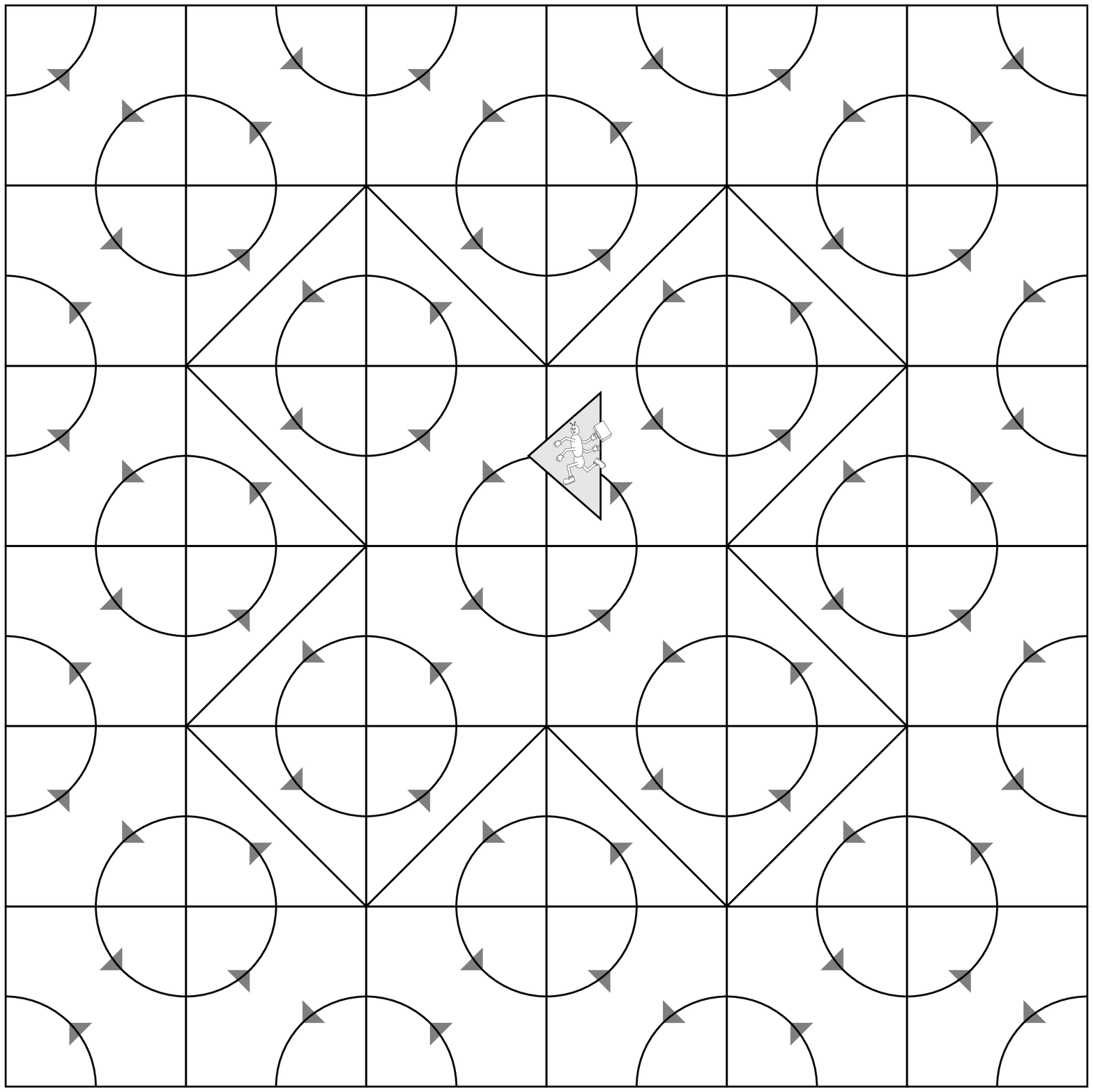,height=2.5in}\hfil
	            \psfig{figure=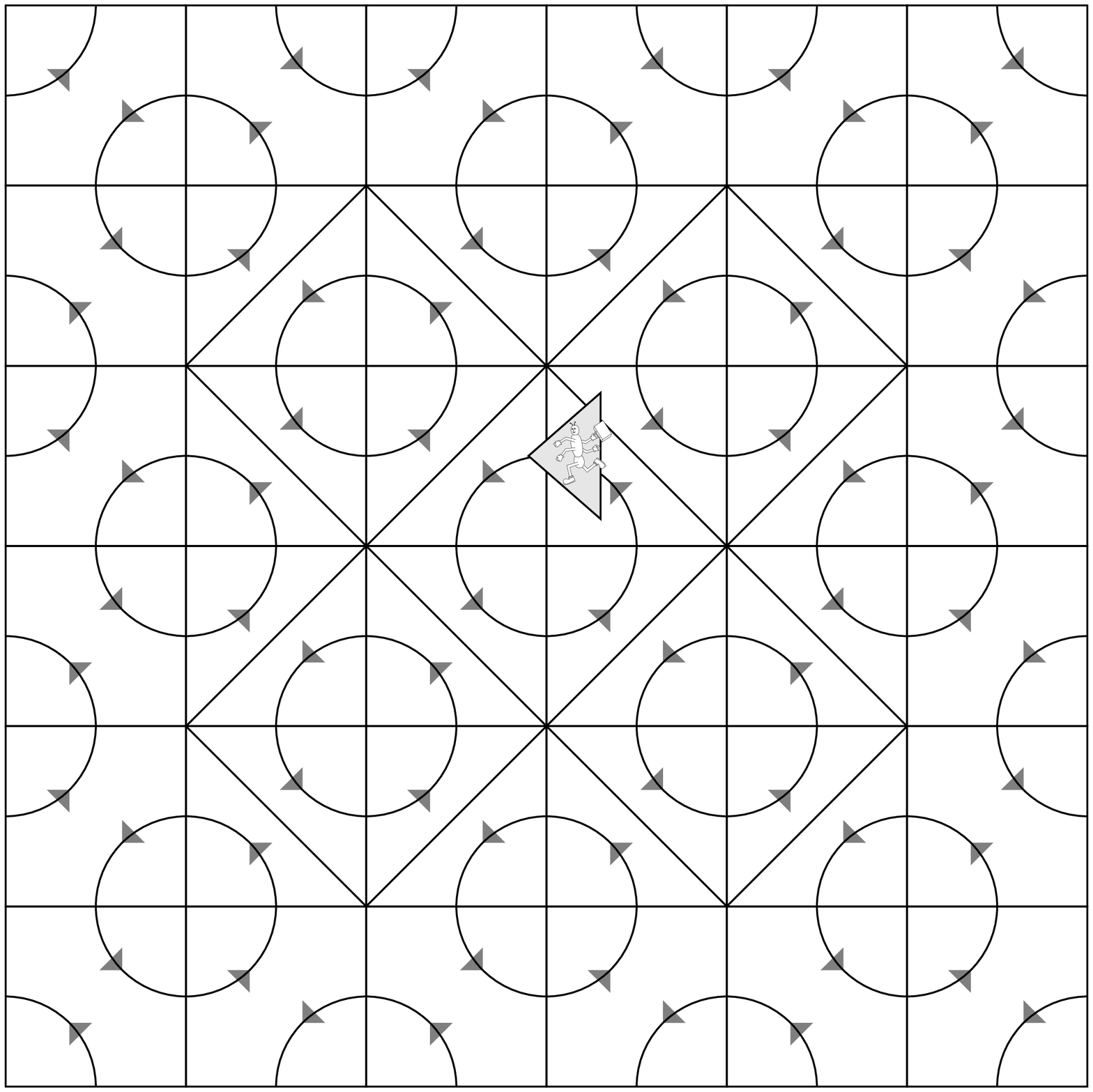,height=2.5in}}
	\caption{The universe of ant 12 at its first four symmetric returns 
                 to home: 4, 8, 28 and 32.} 
\end{figure}

\begin{figure}[htp]
	\centerline{\psfig{figure=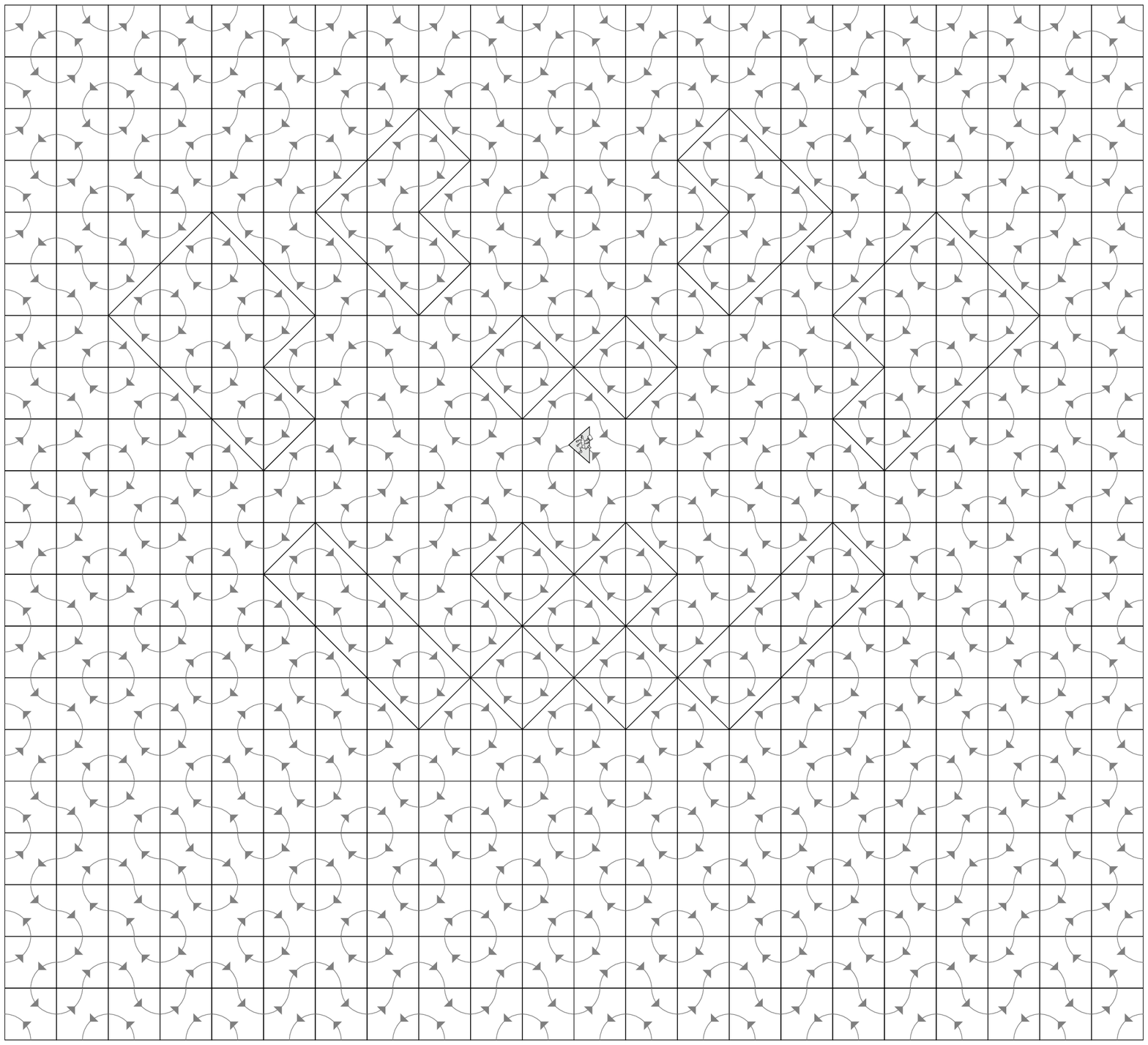,width=\hsize}}
	\caption{ The diagonal graph of ant 48 at time 7016. }
\end{figure}

Why does
the even run-length property imply recurrence of bilateral symmetry?
Here R\"ummler has augmented the picture in a manner that can be
paraphrased as follows.  Let us say that a cell is {\it cold} if its state
is odd, so that it will not change orientation the next time it is
visited (because of the even run-length property), and {\it hot} if its state
is even, so that it may or may not change orientation the next time it
is visited.  To complete the picture we make the convention that for
hot cells we display not only the Truchet tile but also its diagonal.
The {\it diagonals graph} is the graph whose edges consist 
of the diagonals of the hot tiles.  Figures 11 and 12 show the diagonals 
graph of ants 12 and 48 associated with certain instants in time at which 
the ant has returned to its initial location (hereafter called ``home'').  
Note that these graphs can be quite complicated, breaking up into six
components as shown in Figure 12.  However, observe a key fact, 
which we call the {\it even diagonal-degree property}:

\begin{equation}
{\mbox{\it  All vertices in the diagonals graph have even degree (0, 2,
or 4).}}
\end{equation}

\begin{lemma}
Suppose that the ant is at its home position, and that the
state of the universe satisfies (2).  Then the ant will travel along
the principal contour (and return home).
\end{lemma}

{\bf  Proof:} As was remarked earlier, the only thing we have to worry
about is that the principal contour $C$ might visit a cell twice, and
that this cell might change its orientation after the ant's first
visit.  If the twice-visited cell is cold, then the Truchet tile
will not change its orientation after the ant's first visit.  What
about twice-visited cells that are hot?  We will show that such cells
do not exist, as a consequence of (2).  Let $T$ be such a cell, and let
$d$ be the diagonal of $T$ connecting vertices $u$ and $v$. 

\bigskip
{\bf Claim:} If $d$ is deleted then in the resulting graph the
components of $u$ and $v$ are disjoint. 

\bigskip
If we can show this, then the desired contradiction follows, since from (2)
the component of $u$ (or $v$) would contain only one vertex of odd degree,
contradicting the well-known fact (usually associated with Euler) that
a connected graph must always have an even number of vertices of odd
degree.  (This is sometimes referred to as the handshake theorem, in
that it says that the number of people who have shaken hands an odd
number of times is even.)   

\begin{figure}[htp]
    \centerline{\hskip .5in \psfig{figure=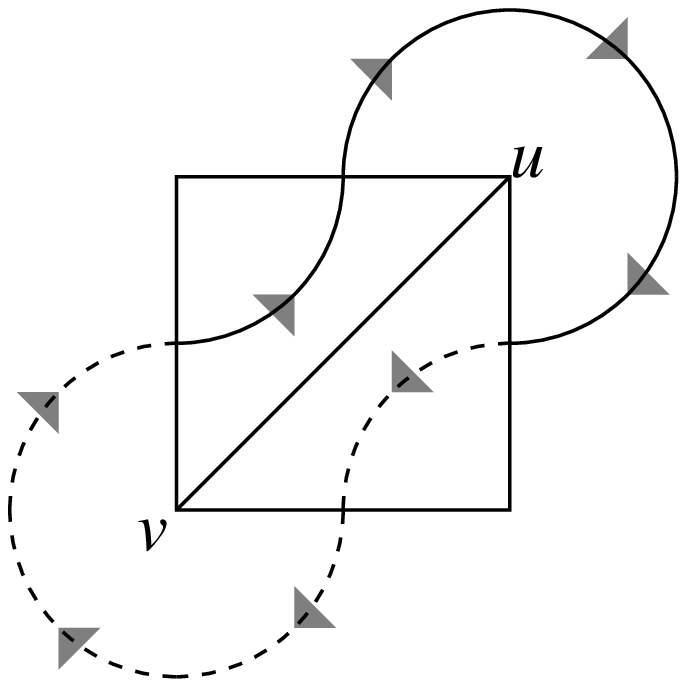,width=1in} \hskip 1.5in
                             \psfig{figure=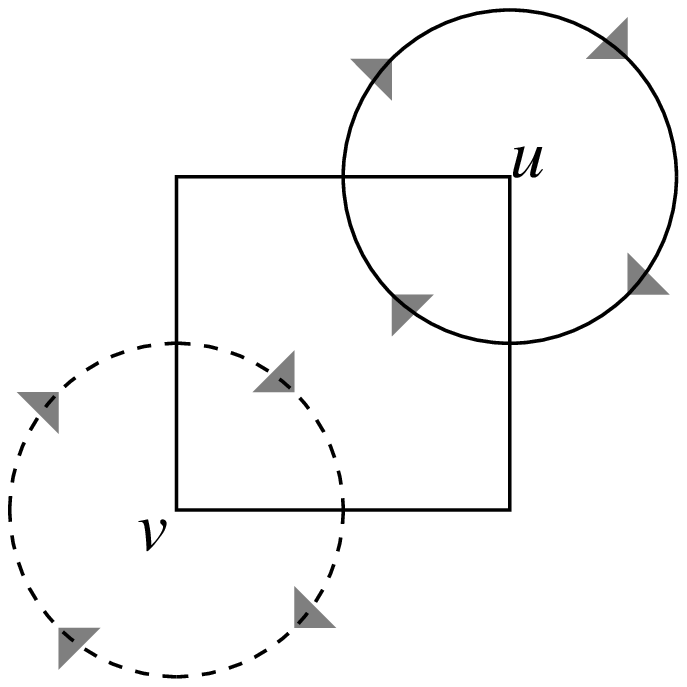,width=1in} }
	\caption{(a) Before \hskip 1.25in (b) After}
\bigskip
    \centerline{\psfig{figure=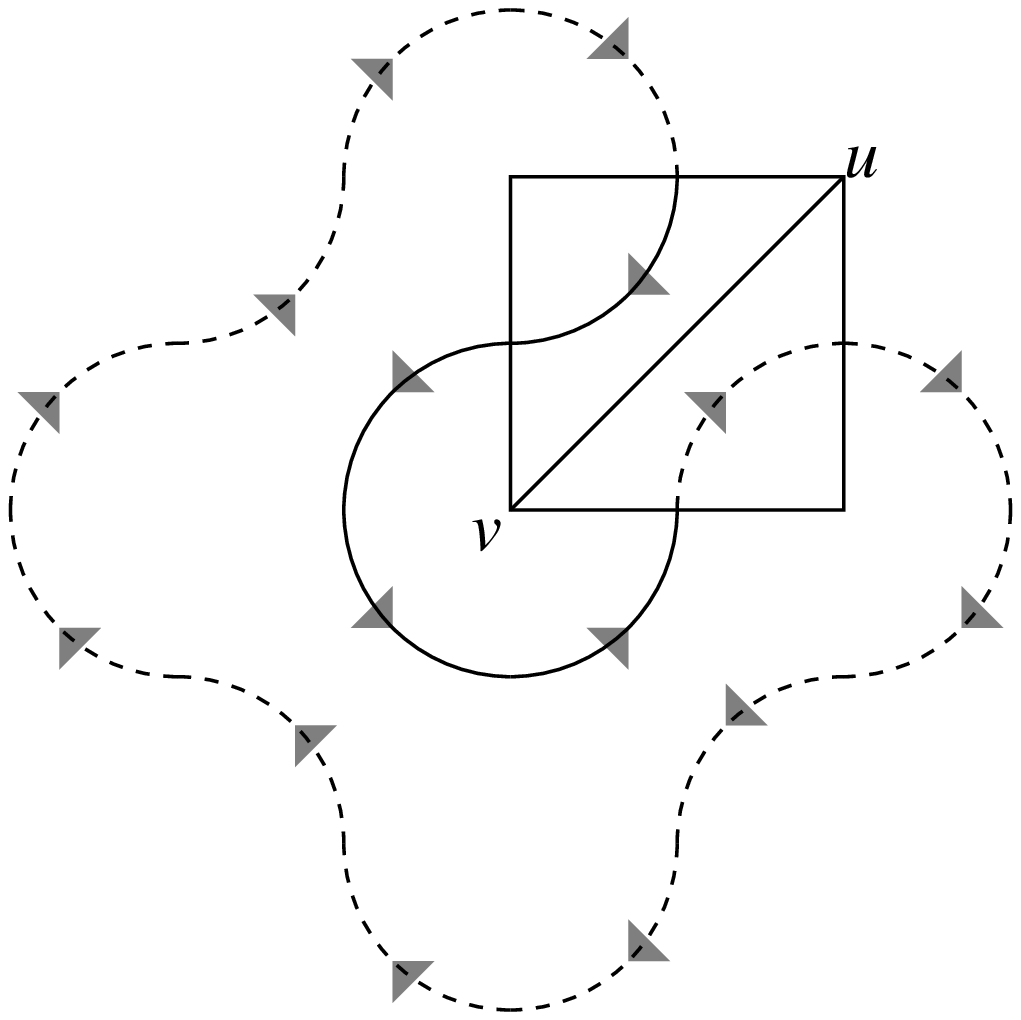,width=1.5in} \hskip 1in
                \psfig{figure=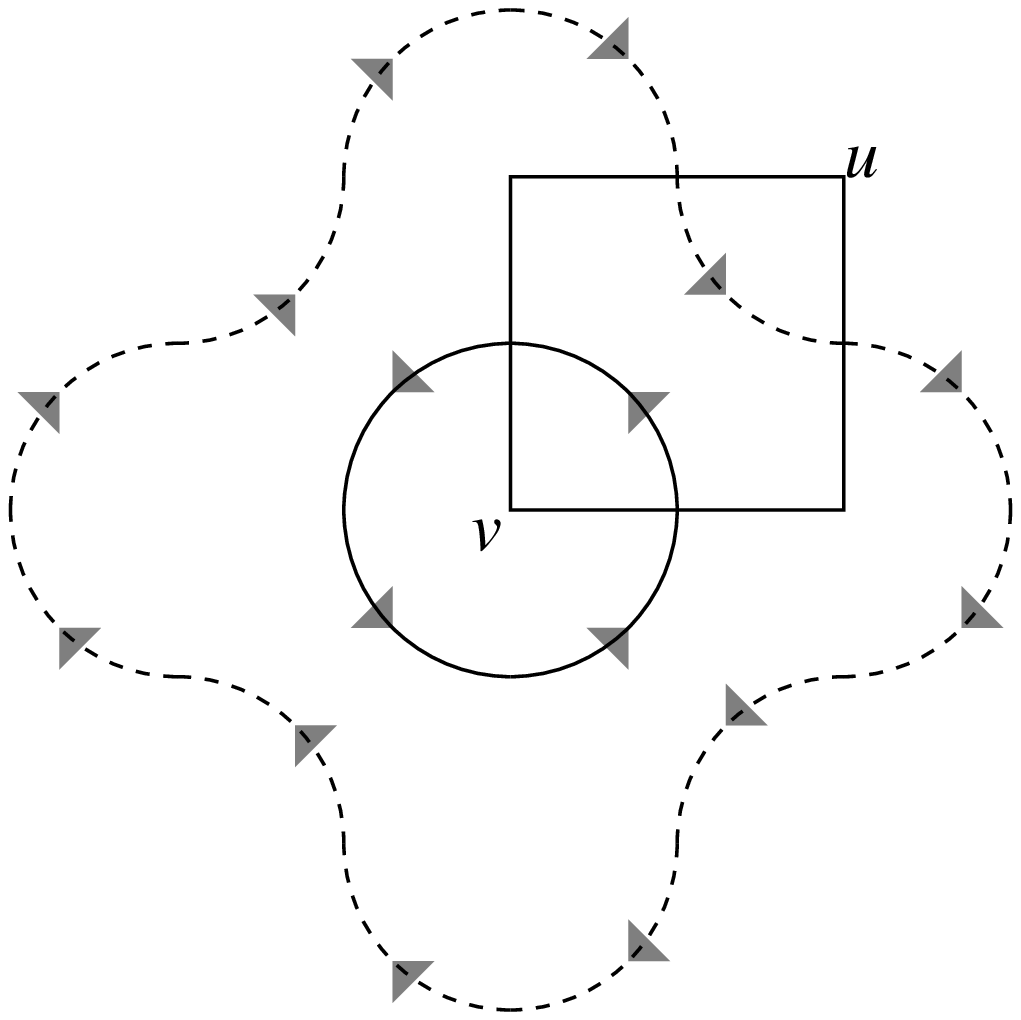,width=1.5in} }
	\caption{(a) Before \hskip 1.25in (b) After }
\end{figure}

It remains to prove the claim.  For this purpose consider the 
twice-visited tile $T$ and its two arcs (quarter circles). 
Without loss of generality, we assume that $T$ is an H-tile.
Now color in red
the arc in $T$ below the diagonal and in addition all succeeding arcs of 
$C$ up to the point where $C$ is about to re-enter $T$.  Color the remaining 
arcs blue.  The dotted arcs in Figures 13(a) and 14(a)
represent the red arcs and the
solid arcs represent the blue.  Now consider the same picture except 
that the diagonal of $T$ has been deleted and $T$ has switched to its other 
orientation, as shown in Figures 13(b) and 14(b); as before, the arc in
$T$ below the diagonal should be colored red and the other arc blue.
One sees at once that $C$ has split into two contours, one all 
red, the other all blue (this much is a purely combinatorial fact and 
has nothing to do with the topology of the plane).  Now the Jordan curve 
theorem tells us that each of these non-intersecting contours has an 
inside and an outside, which can be arranged either as in Figure 13(b)
(the non-nested case) or as in Figure 14(b) (the nested case).  In
either case it is clear that the components of $u$ and $v$ are disjoint,
since both the red contour and blue contour intervene.  Combining this 
with the handshake theorem completes the proof.

\bigskip
{\bf Remark}: While these last observations can be made rigorous without
using the full force of the Jordan curve theorem, there must be some
use of the topology of the plane since the analogue of Lemma 1 is not
true on the torus.

 Finally, we must prove (2).

\begin{lemma}
If (2) holds when the ant is at home it will still hold after
the ant has toured the principal contour and returned home.
\end{lemma}

{\bf Proof}.  Let $v$ be some vertex of the diagonals graph.  The 
{\it neighborhood}, $N(v)$, 
consists of the four cells having $v$ as a vertex.  We will show
that if at some point the Truchet contour enters and then leaves $N(v)$,
the parity of the degree of $v$ is not changed.
   
\bigskip
{\bf Case I}.  The contour meets $N(v)$ 
in only one cell.  Then $v$ does not lie
on a diagonal of that cell.  If the cell is cold the situation remains
the same since cold cells do not switch, and if it is hot then it
becomes cold, hence has no solid diagonal; so, whether it switches or
not, the degree of $v$ is unchanged.

\begin{figure}[htp]
	\centerline{\hfil \psfig{figure=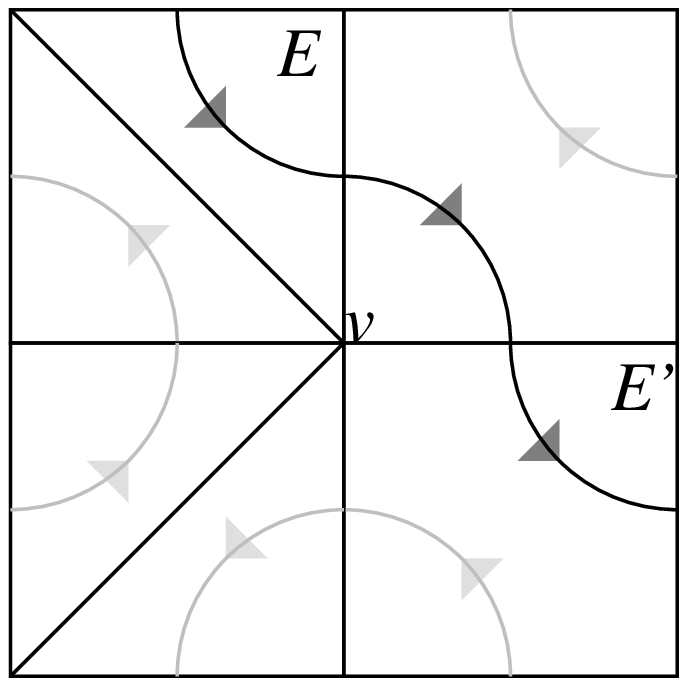,height=1.5in} \hfil
            \psfig{figure=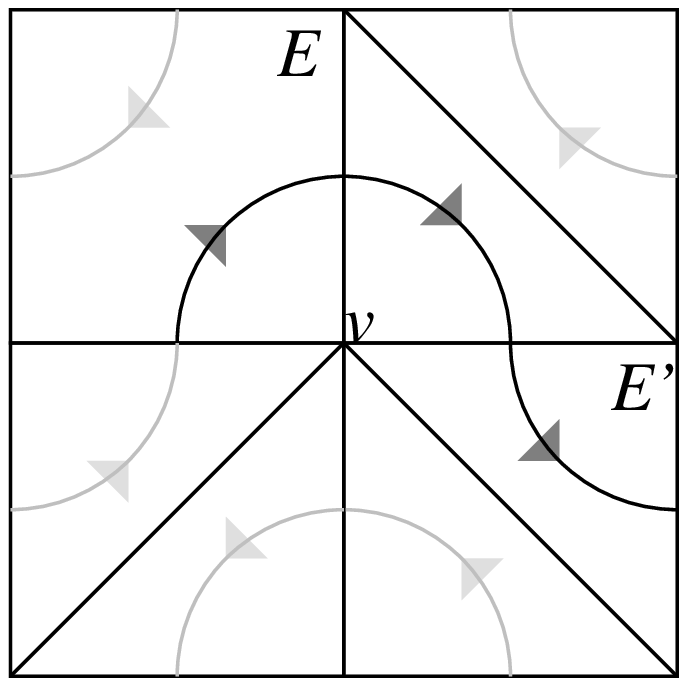,height=1.5in} \hfil}
	\caption{The parity of the number of solid diagonals incident to $v$
	         does not change.}
\end{figure}

\bigskip
{\bf Case II}.  The contour meets more than one cell of $N(v)$.  Then let $E$
be the cell where it enters and $E'$ the cell where it exits.  If $E$ is
hot then its diagonal is incident on $v$ so when it becomes cold this
diagonal disappears (whether $E$ switches or not), and if it is cold it
becomes hot (without switching) so a new solid diagonal will be
incident on $v$.  The exact same argument applies to $E'$, so the net
effect of the two changes is to preserve the parity of the degree of $v$
(see Figure 15).  As for intermediate cells (there may be 1 or 2 of
them), their diagonals do not meet $v$ so the argument is the same as for
Case I: these cells do not contribute to any change in the number of 
solid diagonals incident on $v$.  Note that there is an additional
case, in which the cells $E$ and $E'$ coincide; we leave the analysis
of this case to the reader.

\bigskip
In general, the principal contour may re-enter and re-exit $N(v)$ several
times; but if one looks at the portion of the curve between any 
entrance-point and the corresponding exit-point, the preceding analysis
will apply.  This completes the proof.  

\bigskip
    Combining Lemmas 1 and 2, we can finally see what is happening: If
the state of the universe satisfies the even diagonals-degree
property with the ant at home, then the ant must travel along the 
principal contour, but when it completes this path and returns home,
it restores the even diagonals-degree property, so that it must once 
again travel along the (new) principal contour, and so on, ad infinitum.  
  
   We leave it to the reader to find the simple number theoretic
argument that shows that a number like 57 whose binary expansion has
the even run-length property is divisible by three; this explains
Propp's observation concerning the code-numbers of the ants whose
tracks exhibit recurrent bilateral symmetry.  

\bigskip
{\it Note: For a copy of Propp's program ant.c (an ant-universe 
simulator designed for UNIX machines), send email to {\tt
propp@math.mit.edu}, or see
{\tt http://www.math.sunysb.edu/\~{}scott/ants}. 

\begin{center}
{\bf References}
\end{center}

\begin{enumerate}
\item D. Gale ``The Industrious Ant'', Mathematical Intelligencer, 
vol.~15, no.~2 (1993), pp 54--58.

\item D. Gale and J. Propp ``Further Ant-ics'', Mathematical Intelligencer, 
vol.~16, no.~1 (1994), pp. 37--42.

\item L.A. Bunimovich and S. Troubetzkoy 
``Recurrence properties of Lorentz Lattice Gas Cellular Automata'',
Journal of Statistical Physics, vol.~67 (1992), pp. 289--302.

\end{enumerate}
\end{document}